\documentclass{IEEEtran}

\usepackage{graphicx} 
\graphicspath{{./jpeg-images/}}

\usepackage{amsmath}
\usepackage{caption}
\usepackage{amssymb}
\usepackage{subcaption}
\usepackage{graphics} 
\usepackage{epsfig} 
\usepackage{times} 
\usepackage{xcolor} 

\usepackage{amsmath} 
\usepackage{amssymb,xcolor}  
\usepackage{mathtools,gensymb,amsfonts,amsthm,mathrsfs,comment,bm}
\usepackage[capitalise]{cleveref}

\theoremstyle{plain}
\newtheorem{theorem}{Theorem}
\newtheorem{lemma}[theorem]{Lemma}

\theoremstyle{definition}

\theoremstyle{remark}
\newtheorem{remark}{Remark}

\DeclarePairedDelimiter\abs{\lvert}{\rvert}
\DeclarePairedDelimiter\norm{\lVert}{\rVert}
\makeatletter
\let\oldabs\abs
\def\abs{\@ifstar{\oldabs}{\oldabs*}}
\let\oldnorm\norm
\def\norm{\@ifstar{\oldnorm}{\oldnorm*}}
\makeatother

\title{\LARGE \bf
Two Strings with a Dynamic Interior Mass: A  Feedback Control Design with Guaranteed Exponential Decay}

\author{Zoe Brown and Ahmet \"Ozkan \"Ozer 
\thanks{\small *The second author acknowledges the support of the Fulbright U.S. Scholar Program under the 2024-2025 Research Award to conduct research in France.}
	\thanks{Department of Mathematics, Western Kentucky University (WKU), Bowling Green, KY 42101, USA. 
{\tt\small Email: ozkan.ozer@wku.edu}
	}
}

\begin{document}

\maketitle
\thispagestyle{empty}
\pagestyle{empty}

\begin{abstract}
This paper investigates the exponential stabilization of a coupled two-string system joined by a dynamic interior mass. The combined effect of three feedback mechanisms, boundary damping from tip velocity, higher-order nodal damping from angular velocity, and lower-order nodal damping from mass velocity, is analyzed using a Lyapunov framework. Exponential stability is established unconditionally, without constraints on wave speeds or mass location, improving upon earlier results that lower-order nodal damping, as in {Hansen-Zuazua'95}, or boundary damping alone, as in {Lee-You'89}, does not ensure exponential decay without additional structural conditions. Moreover, the lower-order feedback can be removed without loss of exponential decay when combined with the other two mechanisms, via a compact perturbation argument. These results apply to hybrid systems with interior or tip mass interfaces, including overhead cranes, deep-sea cables, and fluid structure interaction. Theoretical findings are validated through numerical simulations.
\end{abstract}

\begin{IEEEkeywords}
Serially-connected strings, Hybrid PDE–ODE systems, Boundary and interface feedback control, Interior-mass dynamics, Unconditional exponential stability
\end{IEEEkeywords}

\section{Introduction}

{Hybrid PDE–ODE systems arise in many engineering applications, including overhead cranes with flexible cables \cite{coron2007control}, deep sea cables with tip dynamics \cite{WangKrstic2020}, bridge cables with internal or external dampers \cite{di2020cable}, and fluid–structure interaction systems \cite{ammari2017feedback}. These systems often involve interior masses or joints, leading to complex transmission conditions and requiring tailored damping strategies. Recent engineering studies have proposed passive or active mass–spring–damper devices for vibration suppression and energy dissipation in such systems. Applications include overhead transmission lines \cite{bukhari2018vibration}, robotic manipulators \cite{QIU2019623}, underwater cables \cite{quan2020dynamics}, and bridge cables \cite{di2020cable}, where point masses serve as absorbers or actuators under active or semi-active control. Beyond classical devices such as Stockbridge dampers, newer approaches incorporate smart materials like piezoelectric and magnetostrictive elements to provide tunable impedance and feedback at interior or boundary locations \cite{Tsetserukou2007}.}

\begin{figure}[htb!]
\centering
\includegraphics[width=0.45\textwidth]{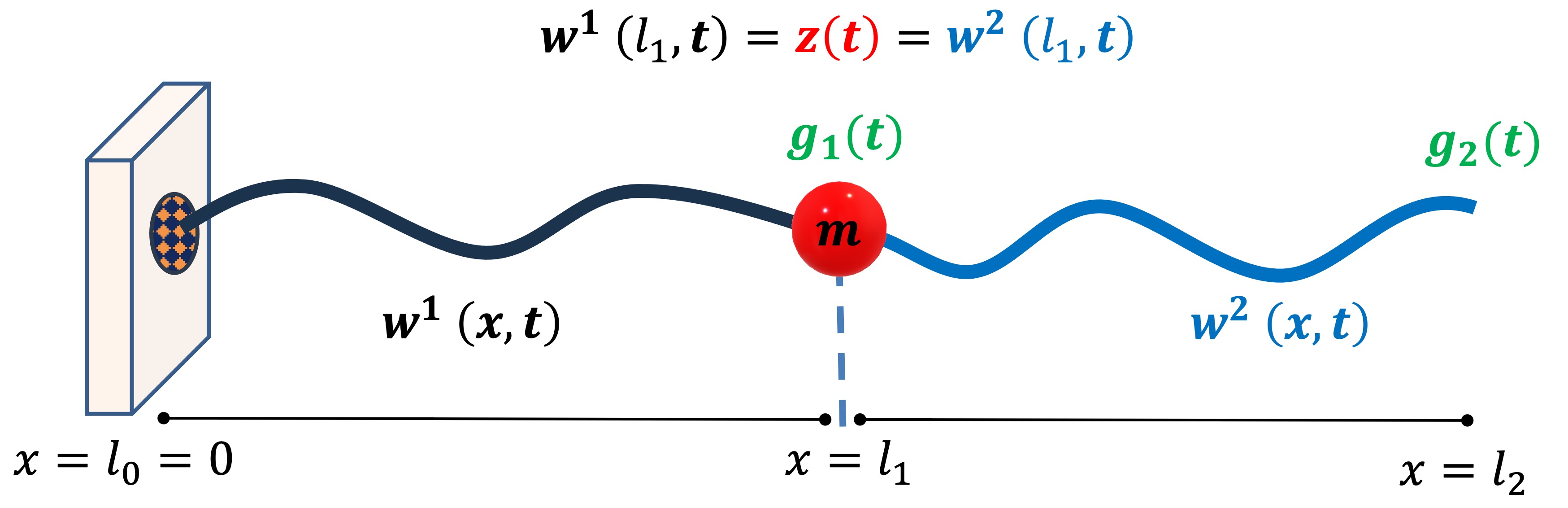}
\caption{\footnotesize Schematic of two coupled strings connected by an interior mass $m$ at $x = l_1$. Control inputs $g_1(t)$ and $g_2(t)$ act at the interface and the right boundary, respectively.}
\label{Intmass}
\end{figure}

This paper investigates one such configuration, focusing on the stabilizing effect of boundary and interface feedback controls. We study a serially-connected system of two nonhomogeneous strings joined by a point mass $m>0$. The dynamics are governed by wave equations modeling transverse vibrations over $(l_0, l_1)$ and $(l_1, l_2)$, with $l_0 = 0$. The interface at $x = l_1$ introduces interior mass dynamics via an inertial term from $m > 0$, coupling the two string segments (see Fig.~\ref{Intmass}). Key parameters include mass densities $\rho_i$ and stiffness coefficients $\alpha_i$ for each segment ($i=1,2$). We prove that the system achieves unconditional exponential stability using a Lyapunov-based approach.

The displacement fields of the two strings are denoted by $w^i(x,t)$ for $i=1,2$. At the coupling interface $x = l_1$, displacement continuity is enforced through a shared dynamic variable $z(t) := w^1(l_1,t) = w^2(l_1,t), $
representing the displacement of the point mass. Let $g_1(t)$ and $g_2(t)$ be the control inputs at the interface $x = l_1$ and the free end $x = l_2$, respectively. The governing equations are given by
\begin{eqnarray}\label{initialSys0}
\left\{\begin{aligned}
    &\rho_i w^i_{tt} -\alpha_i w^i_{xx}(x,t) = 0, && (x,t) \in (l_{i-1}, l_i) \times \mathbb{R}^+,
\end{aligned}\right.\\
\label{initialSys1}
\left\{\begin{aligned}
    &w^1(0,t) = 0,\\
    &w^1(l_1,t) = w^2(l_1,t) = z(t),\\
    &\alpha_1w^1_x(l_1,t) - \alpha_2 w^2_x(l_1,t) + mz_{tt}(t) =g_1(t),\quad \\
    &w^2_x(l_2,t) = g_2(t),
\end{aligned}\qquad~\right.\\
\label{initialSys2}
\left\{\begin{aligned}
    (w^i, w^i_t)(x,0) = (w^i_0, w^i_1)(x),  ~~x \in [l_{i-1},l_i], i=1,2.~
\end{aligned}\right.
\end{eqnarray}
The control inputs are subsequently chosen in the form of feedback laws as defined below
\begin{equation}\label{feedbackLaw}
\left\{
\begin{aligned}
g_1(t) &:= -b_0 \left({ \alpha_1} w^1_{xt}(l_1,t) - { \alpha_2} w^2_{xt}(l_1,t)\right) - b_1 w^1_t(l_1,t),\\
g_2(t) &:= -d_1 w^2_t(l_2,t).
\end{aligned}
\right.
\end{equation}
Specifically, $g_1(t)$ combines higher-order slope-velocity feedback and lower-order interface-velocity feedback at $x = l_1$, while $g_2(t)$ applies lower-order velocity feedback at the free end $x = l_2$. The gains $b_0 > 0$ and $b_1 > 0$ correspond to slope-difference and velocity damping at the interface, respectively, and $d_1 > 0$ denotes boundary damping at the right end. These feedbacks are central to the system’s stabilization, as analyzed in the next sections.

We introduce an auxiliary dynamic variable $\eta(t)$ by
\begin{equation}\label{defeta}
\eta(t) := b_0 \left( {\alpha_1} w^1_x(l_1,t) - {\alpha_2} w^2_x(l_1,t)\right) + m z_t(t),
\end{equation}
By differentiating \eqref{defeta} and using the transmission condition in \eqref{initialSys1}, $\eta(t)$ satisfies the initial-value problem
\begin{equation}
\left\{\begin{aligned}
    &\eta_t(t) = -({ \alpha_1} w^1_x - { \alpha_2} w^2_x)(l_1,t) - b_1 z_t(t),  t \in \mathbb{R}^+,\\
  &\eta(0) = \eta_0:=b_0 \big(\alpha_1 w^1_x- \alpha_2w^2_x\big)(l_1, 0)  + m z_t(0).
\end{aligned}\right.
\end{equation}

\subsection{Literature on the PDE Model and Stability Results}

The stabilization of hyperbolic systems with interior point masses has been widely studied. Early works such as \cite{Chen1987, LeeYou1989} introduced hybrid string models, showing that interior masses prevent strong stability when only velocity-type feedback is applied at the joint ($x=0$), i.e., when $b_0, d_1 \equiv 0, b_1 \neq 0$ in \eqref{feedbackLaw}. If a velocity-type boundary feedback is applied only at the right boundary ($x = \ell_2$), i.e., $b_0,b_1 \equiv 0, d_1 \neq 0$, the system attains strong stability but lacks uniform exponential decay, as eigenvalues cluster near the imaginary axis \cite{H-Z}. {Further analysis in \cite{Littman-Taylor} showed that in such cases, only polynomial decay—typically at rate $t^{-1}$—can be achieved due to the lack of sufficient dissipation.}

{Subsequent studies, including \cite{ammari2017feedback, Castro1, AmmariTucsnak, Avdonin-Edwards, Boughamda}, extended these ideas to more intricate geometries such as networks, trees, and fluid–structure interaction models, often featuring transmission conditions through interior masses or junctions. These works primarily employed lower-order nodal or boundary damping and demonstrated that stabilization depends sensitively on system parameters, interface locations, and controller placement. While strong or polynomial stability can often be achieved under specific configurations, such damping strategies alone generally fail to ensure uniform exponential stabilization.}

{A key advance appeared in \cite{Morgul2, GuoXu2000}, which introduced higher-order (angular velocity) damping for strings with tip masses and established exponential stability via Lyapunov methods and spectral analysis. This marked a shift in control design, showing that higher-order feedback can succeed where lower-order damping fails. Although our approach does not rely on backstepping, we note that backstepping-based boundary control, as introduced in \cite{krstic2007backstepping}, provides a general framework for PDE and PDE–ODE stabilization and is complementary to higher-order strategies. Related work on PDE–ODE–PDE systems, such as deep-sea cables with tip actuators \cite{WangKrstic2020}, further highlights the effectiveness of boundary-based higher-order feedback.

\subsection{Our Motivation and Major Results}

{As shown in \cite{Chen1987, LeeYou1989, H-Z, Littman-Taylor}, neither exponential nor strong stability is guaranteed when using only boundary damping at $x = l_2$ or lower-order damping at the mass interface $x = l_1$. Building on the higher-order feedback approach of \cite{Morgul2}, we develop a new stabilization framework that combines higher- and lower-order damping at the mass with boundary damping at the free end.}

We establish unconditional exponential stability via a Lyapunov-based approach under the assumption that $b_0$, $b_1$, and $d_1$ are positive. {Unlike previous works with $b_0 \equiv 0$ in \eqref{feedbackLaw}, our model includes higher-order angular velocity feedback at the mass, which plays a decisive role in stabilization.} Remarkably, we show that the lower-order damping $b_1$ can be omitted without loss of exponential decay by treating the reduced system as a compact perturbation. {This highlights the spectral dominance of higher-order damping and enables a more efficient design.}

{The proposed strategy is robust concerning wave speeds, damping strengths, and mass location, and does not rely on structural compatibility conditions, an improvement over earlier results.}

The remainder of the paper is organized as follows. Section~\ref{well} establishes well-posedness via a first-order reformulation. Section~\ref{expstab} presents the main exponential stability results. Section~\ref{numerics} introduces a finite difference scheme preserving the decay rate {and} provides numerical simulations validating the theory.

\section{Well-Posedness}
\label{well}
To analyze the stability of the system \eqref{initialSys0}–\eqref{initialSys2}, we define the total energy $E(t)$ as
\begin{equation} \label{energy}
\begin{aligned}
E(t) = \, &\frac{1}{2} \sum\limits_{j=1}^2 \int_{l_{j-1}}^{l_j} \left[ \rho_j (w^j_t(x,t))^2 + \alpha_j (w^j_x(x,t))^2 \right] \, dx \\
&+ \frac{1}{2} \frac{(\eta(t))^2}{m+b_0b_1}.
\end{aligned}
\end{equation}

\begin{lemma}\label{lana}For all $t\ge 0,$ the energy is dissipative,
\begin{equation} \label{dissp}
\begin{aligned} 
\frac{dE(t)}{dt} &=  - \frac{b_0}{m+b_0 b_1}\left[\alpha_1 w^1_x - \alpha_2 w^2_x \right]^2(l_1,t) \\
&\qquad - \frac{mb_1}{m+b_0 b_1} (z_t(t))^2-d^1(w^2_t)^2(l_2,t). \end{aligned} \end{equation}
	\end{lemma}
\begin{proof}
We differentiate \eqref{energy} along with  \eqref{initialSys0}–\eqref{defeta}:
\begin{eqnarray*}
\begin{array}{ll}
\frac{dE(t)}{dt}= \sum\limits_{j=1}^2 \int_{l_{j-1}}^{l_j} \rho_j (w^j_t(x,t)) (w^j_{tt}(x,t)) \, dx \\
\quad  + \sum\limits_{j=1}^2 \int_{l_{j-1}}^{l_j} { \alpha}_j (w^j_x(x,t)) (w^j_{xt}(x,t)) \, dx + \frac{\eta(t) \eta_t(t)}{m+b_0b_1}\\
=\sum\limits_{j=1}^2 \left. { \alpha_j}  (w^j_t(x,t)) (w^j_{x}(x,t)) \right|_{l_{j-1}}^{l_j}+ \frac{\eta(t) \eta_t(t)}{m+b_0b_1}\\
= -d_1 (w^2_t)^2(l_2,t)- { \alpha_2} (w^2_t w^2_{x})(l_1,t)+ { \alpha_1} (w^1_t w^1_{x})(l_1,t)\\
\quad+ \frac{(b_0 ({ \alpha_1} w^1_x - { \alpha_2} w^2_x)(l_1, t) + mz_t(t)) (-({ \alpha_1} w^1_x - { \alpha_2} w^2_x)(l_1, t) - b_1 z_t(t))}{m+b_0 b_1}.
              \end{array}
\end{eqnarray*}

Using $z(t) = w^1(l_1,t) = w^2(l_1,t)$ and simplifying yields the dissipation identity \eqref{dissp}.
\end{proof}

Let  $H^1_L(0,l_1)=\left\{f\in H^1(0,l_1)~ : ~ f(l_1)=0\right\}.$ Now, define the energy space $\mathcal{H}$ by
$
\mathcal{H}=H^1_L(0, l_1)\times L^2(0,l_1)\times H^1(l_1, l_2)\times L^2(l_1, l_2)\times \mathbb{C}
$ equipped with the  inner product
\begin{equation*}
\begin{array}{ll}
\displaystyle
\langle U, \widetilde{U}\rangle _{\mathcal{H}}=&\sum\limits_{j=1}^2\int_{l_{j-1}}^{l_j} \left[ \rho_j (p^j \bar{\tilde{p}}^j)+ \alpha_j (u^j_x \bar{\tilde{u}}^j_x)\right](x,t) dx\\
&+ \frac{\theta \overline{\widetilde{\theta}}}{m+b_0b_1}
\end{array}
\end{equation*}
for $U=(u^1,p^1,u^2,p^2,\theta)$,  $\widetilde{U}=(\tilde{u}^1,\tilde{p}^1, \tilde{u}^2,\tilde{p}^2,\tilde{\theta})\in \mathcal{H}$. 
 The norm induced by the inner product is $\|U\|^2_{\mathcal{H}}:=\langle U, U\rangle _{\mathcal{H}}.$
Introducing $w^1_t = v^1$ and $w^2_t = v^2$, the system \eqref{initialSys0}–\eqref{initialSys2} can be recast as the first-order initial value problem
\begin{eqnarray}\label{IVP}
\frac{d}{dt}{U}(t) = \mathcal{\mathcal{ A}}{U}(t),\quad {U}(0) = {U}_{0},\quad
\forall\, t > 0,
\end{eqnarray}
where $U = (w^1, v^1, w^2, v^2, \eta)^{\top}$ and $U_0 = (w^1_0, w^1_1, w^2_0, w^2_1, \eta_0)^{\top}$. The operator $\mathcal{A} : \mathcal{D}(\mathcal{A}) \subset \mathcal{H} \rightarrow \mathcal{H}$ is defined by
\begin{eqnarray}
\nonumber \begin{array}{ll}
&\mathcal{A}
[w^1, v^1, w^2, v^2, \eta]^\top=[v^1, \dfrac{\alpha_1}{\rho_1} w^1_{xx}, \\
 &\qquad\qquad v^2, \dfrac{\alpha_2}{\rho_2} w^2_{xx}, -(w^1_x - w^2_x)(l_1) - b_1 v^1(l_1)]^\top.
\end{array}
\end{eqnarray}
with domain
\begin{equation*}
\begin{array}{ll}
\mathcal{D}({\mathcal{ A}})=\left\{
U\in {\mathcal{ H}}:  w^1\in \left(H_{L}^{1}\cap H^{2}\right)(0,l_1), \right.\\
w^2\in H^2(l_1,l_2), v^1\in H_{L}^{1}(0,l_1), v^2\in H^{1}(l_1,l_2), \\
\alpha_2w^2_x(l_2)=-d_1v^2(l_2), \\
\left. b_0 (\alpha_1w^1_x(l_1) - \alpha_2w^2_x(l_1)) + m v^1(l_1)=\eta\right\}.
\end{array}
\end{equation*}
If $(w^1,v^1,w^2,v^2,\eta)$ is a sufficiently regular solution of \eqref{initialSys0}–\eqref{initialSys2}, then it satisfies the abstract Cauchy problem
\begin{equation}\label{CP}
U_t = \mathcal{A} U, \quad U(0) = U_0,
\end{equation}
in the Hilbert space $\mathcal{H}$.
where $U = (w^1,v^1,w^2, v^2,\eta)$ and $U_0 = (w^1_0,w^1_1,w^2_0, w^2_1,\eta_0)$.  
A direct computation yields 
\begin{equation*}\label{Au-u}
\begin{array}{ll}
    {\mathrm Re} (\mathcal{A} U, U)_{\mathcal{H}}  =&  - \frac{b_0}{m+b_0 b_1}\left|(\alpha_1 w^1_x - \alpha_2 w^2_x)(l_1) \right|^2 \\
   & - \frac{mb_1}{m+b_0 b_1} |v^1(l_1)|^2-d^1|v^2(l_2)|^2,
    \end{array}
\end{equation*}
showing that $\mathcal{A}$ is dissipative. To prove $m$-dissipativity, let $F = (f_1, f_2, f_3, f_4, f_5) \in \mathcal{H}$ be given. By the Lax–Milgram Theorem, there exists a unique $U \in \mathcal{D}(\mathcal{A})$ solving $-\mathcal{A}U = F$. Thus, $\mathcal{A}$ is $m$-dissipative and $0 \in \rho(\mathcal{A})$.

By the Lumer–Phillips Theorem, $\mathcal{A}$ generates a $C_0$-semigroup of contractions $(e^{t \mathcal{A}})_{t \ge 0}$ on $\mathcal{H}$, and the solution to \eqref{CP} is given by
$(e^{t\mathcal{A}})_{t \geq 0}$ establishing the well-posedness of the system.
\begin{theorem}
Letting $U_0\in \mathcal{H}$, the system \eqref{CP} admits a unique weak solution $U$ satisfying $U \in C^0 (\mathbb{R}_+,\mathcal{H})$.\\
Moreover, if $U_0\in D(\mathcal{A})$, the system \eqref{CP} admits a unique strong solution $U$ satisfying $U \in C^1 (\mathbb{R}_+,\mathcal{H}) \cap C^0 (\mathbb{R}_+,D(\mathcal{A}))  $.
\end{theorem}

\section{Exponential Stability Results}

\label{expstab}
\noindent For $\epsilon_1, \epsilon_2 > 0$, we define the Lyapunov function $L(t)$ as
\begin{equation}\label{lyapu1}
\begin{aligned}
L(t) &:= E(t) + \epsilon_1 I_1(t) + \epsilon_2 I_2(t), \\
V(t) &:= t L(t) + P_1(t) + P_2(t),
\end{aligned}
\end{equation}
where the auxiliary functionals $I_j(t)$ and $P_j(t)$, for $j = 1, 2$, are defined by
\begin{equation*}
\begin{aligned}
I_j(t) &:= 2 \int_{l_{j-1}}^{l_j} \rho_j w^j_t(x,t) w^j_x(x,t) \, dx, \\
P_j(t) &:= 3 \int_{l_{j-1}}^{l_j} (x - l_{j-1}) \rho_j w^j_x(x,t) w^j_t(x,t) \, dx.
\end{aligned}
\end{equation*}

\begin{lemma}\label{equiv2}
Let $C := \max\limits_{j=1,2} \sqrt{\frac{\rho_j}{\alpha_j}}$ and $\epsilon := \max\limits_{j=1,2} \epsilon_j$. Then, for any $0 < \epsilon < \frac{1}{2C}$, $L(t)$ defined in \eqref{lyapu1} satisfies 
\begin{eqnarray}
\begin{array}{ll}
(1-2\epsilon C)E(t)\le L(t)\le (1+2\epsilon C) E(t).
\end{array}
\end{eqnarray}
\end{lemma}
\begin{proof}
For each $j = 1, 2$, we apply Young's inequality to estimate the cross term in $I_j(t)$:
\begin{equation*}
\left|I_j(t)\right| \le \sqrt{\frac{\rho_j}{\alpha_j}} \int_{l_{j-1}}^{l_j} \left[ \rho_j (w^j_t)^2(x,t) + \alpha_j (w^j_x)^2(x,t) \right] dx.
\end{equation*}
Summing over $j=1,2$ and recalling the definition of $L(t)$, the desired estimate follows immediately.
\end{proof}

\begin{lemma}\label{equiv3}
Let $C_1 := 3 \max\limits_{j=1,2} \left( (l_j - l_{j-1}) \sqrt{\frac{\rho_j}{\alpha_j}} \right) > 0$, and define $C_2 := \frac{C_1}{1 - 2\epsilon C}$. Then for any $0 < \epsilon < \frac{1}{2C}$, the functional $V(t)$ defined in \eqref{lyapu1} satisfies 
\begin{eqnarray*}
\begin{array}{ll}
(t-C_2)L(t)\le  V(t)\le (t+C_2) L(t),\qquad \forall t\ge 0.\\
\end{array}
\end{eqnarray*}
\end{lemma}
\begin{proof}
By Young's inequality, we estimate each $P_j(t)$ term as follows
\begin{eqnarray*}
\begin{array}{ll}
\left|P_1(t)\right|\le \frac{3 l_1}{2} \sqrt{\frac{\rho_1}{\alpha_1}} \int_{l_0}^{l_1} \left(\rho_1 (w^1_t)^2 +\alpha_1 (w^1_x)^2\right) dx,\\
\left|P_2(t)\right|\le \frac{3 (l_2-l_1)}{2} \sqrt{\frac{\rho_2}{\alpha_2}} \sum\limits_{j=1}^2\int_{l_1}^{l_1} \left(\rho_2(w^2_t)^2 +\alpha_2 (w^2_x)^2\right) dx.
\end{array}
\end{eqnarray*}
Summing the two estimates yields
$
|P_1(t)+P_2(t)|\le C_1 E(t)\le  \frac{C_1}{1-2\epsilon C} L(t),$ where the second inequality uses Lemma~\ref{equiv2}, completing the proof of Lemma~\ref{equiv3}.
\end{proof}
\begin{lemma} \label{dLdt}
The Lyapunov functional $L(t)$ defined in \eqref{lyapu1} satisfies $\frac{dL(t)}{dt} \le 0$ provided that the feedback gains $b_0, b_1, d_1 > 0$ and the positive parameters $\delta$, $\epsilon_1$, and $\epsilon_2$ satisfy the following conditions
\begin{equation}\label{asscon1}
\begin{aligned}
& \epsilon_2 \le \min\left( \frac{1}{2C}, \frac{\alpha_2 d_1}{d_1^2 + \alpha_2 \rho_2} \right), \\
&\epsilon_1 \le \min\left( \frac{1}{2C}, \frac{\epsilon_2 b_0 \alpha_1}{b_0 \alpha_2 + \epsilon_2 (m + b_0 b_1)}, \frac{m b_1}{\rho_1 (m + b_0 b_1)} \right), \\
&\frac{b_0 \alpha_1}{b_0 \alpha_2 + (m + b_0 b_1)\epsilon_2} < \delta < \frac{b_0 \alpha_1 - (m + b_0 b_1)\epsilon_1}{b_0 \alpha_2}.
\end{aligned}
\end{equation}
\end{lemma}

\begin{proof} For $j = 1,2$, differentiating $I_j(t)$ along the solutions of \eqref{initialSys0}–\eqref{initialSys2} yields
\begin{align}
\frac{dI_j(t)}{dt} &= 2 \int_{l_{j-1}}^{l_j} \alpha_j w^j_{xx} w^j_x \, dx + 2 \int_{l_{j-1}}^{l_j} \rho_j w^j_t w^j_{xt} \, dx \notag \\
&= \left. \left[\alpha_j (w^j_x)^2 + \rho_j (w^j_t)^2\right] \right|_{l_{j-1}}^{l_j}. \label{eq:dIjt}
\end{align}
Summing over $j$ and simplifying
\begin{align}
&\sum_{j=1}^2 \epsilon_j \frac{dI_j(t)}{dt} 
\le \epsilon_2 \left( \frac{d_1^2}{\alpha_2} + \rho_2 \right)(w^2_t)^2(l_2) \notag \\
&\quad + (\alpha_1 \epsilon_1 (w^1_x)^2 - \alpha_2 \epsilon_2 (w^2_x)^2)(l_1) + \epsilon_1 \rho_1 (z_t)^2. \label{eq:sumI}
\end{align}
	Using this and Lemma~\ref{lana}, the derivative of $L(t)$ becomes
\begin{align}
\frac{dL(t)}{dt} 
=& - \frac{b_0}{m+b_0 b_1} \left[\alpha_1 w^1_x - \alpha_2 w^2_x\right]^2(l_1) \notag \\
& - \left[\frac{mb_1}{m+b_0 b_1} - \epsilon_1 \rho_1\right] (z_t)^2 \notag \\
&- \left[d_1 - \epsilon_2 \left(\frac{d_1^2}{\alpha_2} + \rho_2\right)\right](w^2_t)^2(l_2) \notag \\
& +  (\alpha_1 \epsilon_1 (w^1_x)^2 - \alpha_2 \epsilon_2 (w^2_x)^2)(l_1). \label{eq:lyapunov_termwise}
\end{align}
We now expand the interface term via the identity
$
\left[\alpha_1 w^1_x - \alpha_2 w^2_x\right]^2 \ge (\alpha_1)^2 (w^1_x)^2 + (\alpha_2)^2 (w^2_x)^2
 - \delta \alpha_1 \alpha_2 (w^1_x)^2 - \frac{1}{\delta} \alpha_1 \alpha_2 (w^2_x)^2.$
Inserting this into \eqref{eq:lyapunov_termwise} gives
\begin{align}
&\frac{dL(t)}{dt} \le - \left[\frac{b_0 (\alpha_1)^2}{m + b_0 b_1} - \alpha_1 \epsilon_1 - \frac{\delta b_0 \alpha_1 \alpha_2}{m + b_0 b_1} \right](w^1_x)^2(l_1) \notag \\
&- \left[\frac{b_0 (\alpha_2)^2}{m + b_0 b_1} + \alpha_2 \epsilon_2 - \frac{b_0 \alpha_1 \alpha_2}{\delta (m + b_0 b_1)} \right](w^2_x)^2(l_1) \notag \\
& - \left[\frac{mb_1}{m + b_0 b_1} - \epsilon_1 \rho_1\right] (z_t)^2 \notag  - \left[d_1 - \epsilon_2 \frac{d_1^2+\alpha_2\rho_2}{\alpha_2}\right](w^2_t)^2(l_2).
\end{align}
Each bracket is nonnegative under \eqref{asscon1}, so $\frac{dL(t)}{dt} \le 0$.
	\end{proof}
	\begin{theorem}\label{majres0}
Let the parameters $\delta, \epsilon_1, \epsilon_2>0$ satisfy
\begin{equation}\label{asscon2}
\left\{
\begin{aligned}
&\epsilon_2 \le \min\left( \frac{1}{2C}, \sqrt{\frac{\alpha_2}{\rho_2}}, \frac{\alpha_2 d_1}{d_1^2 + \alpha_2 \rho_2} \right), \\
&\epsilon_1 \le \min\left( \frac{1}{2C}, \sqrt{\frac{\alpha_1}{\rho_1}}, \right.\\
&\left. \qquad \qquad \frac{\epsilon_2 b_0 \alpha_1}{b_0 \alpha_2 + \epsilon_2 (m + b_0 b_1)}, \frac{m b_1}{\rho_1 (m + b_0 b_1)} \right), \\
&\frac{b_0 \alpha_1}{b_0 \alpha_2 + (m + b_0 b_1)\epsilon_2} < \delta < \frac{b_0 \alpha_1 - (m + b_0 b_1) \epsilon_1}{b_0 \alpha_2}.
\end{aligned}
\right.
\end{equation}
Then, there exists a time $T > \max(T_1, T_2, T_3, T_4)$ such that $\frac{dV(t)}{dt} \le 0$ for all $t \ge T$, where\\
$T_1:= \frac{ \tfrac{3}{2} \alpha_1 l_1 (m + b_0 b_1) + 2 b_0^2 \alpha_1^2 }{ b_0 \alpha_1^2 - \alpha_1 \epsilon_1 (m + b_0 b_1) - \delta b_0 \alpha_1 \alpha_2 }, $\\
$T_2 := \frac{ \tfrac{3}{2} \alpha_2 (l_2 - l_1) (m + b_0 b_1) + 2 b_0^2 \alpha_2^2 }{ b_0 \alpha_2^2 + \alpha_2 \epsilon_2 (m + b_0 b_1) - \frac{b_0 \alpha_1 \alpha_2}{\delta} },$ \\
$T_3 := \frac{ \tfrac{3}{2} \rho_1 l_1 (m + b_0 b_1) + m^2 }{ m b_1 - \epsilon_1 \rho_1 (m + b_0 b_1) }, 
T_4 := \frac{ \tfrac{3}{2} (l_2 - l_1) \rho_2 \alpha_2 }{ \alpha_2 d_1 - \epsilon_2 (d_1^2 + \alpha_2 \rho_2) }.$
\end{theorem}

\begin{proof}
For each $j=1,2$, we compute $\frac{d}{dt}P_j(t)$ by differentiating and integrating by parts
\begin{equation*}
\begin{aligned}
&\frac{dP_j(t)}{dt} 
= \tfrac{3}{2}(l_j - l_{j-1}) \left[ \rho_j (w^j_t(l_j,t))^2 + \alpha_j (w^j_x(l_j,t))^2 \right] \\
&\quad - \tfrac{3}{2} \int_{l_{j-1}}^{l_j} \left[ \rho_j (w^j_t)^2 + \alpha_j (w^j_x)^2 \right] dx.
\end{aligned}
\end{equation*}
Combining this with the Lyapunov structure, and using $\epsilon_1 < \sqrt{\frac{\alpha_1}{\rho_1}}$ and $\epsilon_2 < \sqrt{\frac{\alpha_2}{\rho_2}}$, we estimate
\begin{equation*}
\begin{aligned}
&\frac{dV(t)}{dt} 
= \tfrac{1}{2} \sum_{j=1}^2 \int_{l_{j-1}}^{l_j} \left[ \rho_j (w^j_t)^2 + \alpha_j (w^j_x)^2 \right] dx  \\
&\quad + \tfrac{1}{2(m + b_0 b_1)} \eta(t)^2 + \epsilon_1 I_1(t) + \epsilon_2 I_2(t) \\
&\quad - \left[t \left(\tfrac{b_0 \alpha_1^2}{m + b_0 b_1} - \alpha_1 \epsilon_1 - \tfrac{\delta b_0 \alpha_1 \alpha_2}{m + b_0 b_1} \right) - \tfrac{3}{2} \alpha_1 l_1 \right] (w^1_x(l_1,t))^2 \\
&\quad - \left[t \left(\tfrac{b_0 \alpha_2^2}{m + b_0 b_1} + \alpha_2 \epsilon_2 - \tfrac{b_0 \alpha_1 \alpha_2}{\delta(m + b_0 b_1)} \right)\right.\\
&\left. \qquad\qquad- \tfrac{3}{2} \alpha_2 (l_2 - l_1) \right] (w^2_x(l_1,t))^2 \\
&\quad - \left[ t \left( \tfrac{m b_1}{m + b_0 b_1} - \epsilon_1 \rho_1 \right) - \tfrac{3}{2} \rho_1 l_1 \right] (z_t(t))^2 \\
&\quad - \left[ t \left( d_1 - \epsilon_2 \left( \tfrac{d_1^2}{\alpha_2} + \rho_2 \right) \right) - \tfrac{3}{2} \rho_2 (l_2 - l_1) \right] (w^2_t(l_2,t))^2.
\end{aligned}
\end{equation*}
To ensure nonpositivity, we introduce additional buffer terms in the inequalities
\begin{equation*}
\begin{aligned}
\frac{dV(t)}{dt} &\le - \left[ tA_1 - \tfrac{3}{2} \alpha_1 l_1 - \tfrac{2 b_0^2 \alpha_1^2}{m + b_0 b_1} \right] (w^1_x(l_1,t))^2 \\
&\quad - \left[ tA_2 - \tfrac{3}{2} \alpha_2 (l_2 - l_1) - \tfrac{2 b_0^2 \alpha_2^2}{m + b_0 b_1} \right] (w^2_x(l_1,t))^2 \\
&\quad - \left[ tA_3 - \tfrac{3}{2} \rho_1 l_1 - \tfrac{m^2}{m + b_0 b_1} \right] (z_t(t))^2 \\
&\quad - \left[ tA_4 - \tfrac{3}{2} \rho_2 (l_2 - l_1) \right] (w^2_t(l_2,t))^2,
\end{aligned}
\end{equation*}
where $A_1$, $A_2$, $A_3$, and $A_4$ are the coefficients inside the time-multiplied terms. 
Under the assumption $t > \max(T_1, T_2, T_3, T_4)$ from Theorem~\ref{majres0}, each bracketed term is positive, hence $\frac{dV(t)}{dt} \le 0$ for all $t \ge T$.
\end{proof}

\begin{theorem} \label{majres}
Assume that the positive constants $\delta, \epsilon_1, \epsilon_2$ satisfy the conditions stated in Theorem~\ref{majres0}. Then, there exist constants $M, \sigma > 0$ such that the total energy \eqref{energy} satisfies the exponential decay estimate:
\begin{equation} \label{finalresult}
E(t) \le M e^{-\sigma t} E(0), \quad \forall t > 0.
\end{equation}
\end{theorem}
\begin{proof}
By Lemma~\ref{equiv3}, for all $t > T$, we have
$
(t - C_2) L(t) \le V(t) \le V(T) \le (T + C_2) L(0),
$
which implies
$
L(t) \le \frac{T + C_2}{t - C_2} L(0).
$
Combining this with Lemma~\ref{equiv2}, we deduce:
\begin{equation} \label{exp_decay_bound}
\begin{split}
E(t) &\le \frac{L(t)}{1 - 2\epsilon C}
\le \frac{1 + 2\epsilon C}{1 - 2\epsilon C}  \frac{T + C_2}{t - C_2} E(0).
\end{split}
\end{equation}
Therefore, there exists a time $T^* > T$ and a constant $0 < \zeta < 1$ such that $E(t) \le \zeta E(0)$ for all $t > T^*$. The exponential decay estimate \eqref{finalresult} then follows by standard semigroup arguments.
\end{proof}
{Before the main result, we note that the main challenge is the coupling at the interior mass. The interface variable $\eta(t)$ allows us to capture this in the Lyapunov analysis. Exponential decay is shown even without the lower-order damping $b_1$, via a compactness argument.}

\begin{theorem} \label{majres2} Consider the system \eqref{initialSys0}–\eqref{initialSys2} with $b_0, d_1 > 0$ and $b_1 = 0$. Then, there exist constants $\tilde M, \tilde \sigma > 0$ such that the total energy \eqref{energy} satisfies \begin{equation} \label{finalresult2} E(t) \le \tilde M e^{-\tilde \sigma t} E(0), \quad \forall t > 0. \end{equation} \end{theorem}
\begin{proof}[Sketch of the Proof]
Theorem~\ref{majres} establishes exponential stability when $b_1 > 0$. Removing $b_1$ modifies the generator $\mathcal{A}$ by a bounded finite-rank perturbation, as the damping $b_1 w^1_t(\ell_1,t)$ is localized. Hence, $\left.\mathcal{A}\right|_{b_1=0}$ is a compact perturbation of $\mathcal{A}$. By \cite[Theorem 3.2]{Triggiani1975}, exponential stability carries over, yielding \eqref{finalresult2}.
\end{proof}

\begin{figure}[htb]
    \centering
    \begin{minipage}{0.4\linewidth}
        \centering
        \small \textbf{(a)} $b_0 = 0$, $b_1 = d_1 = 1$\\[0.5em]
        \includegraphics[width=\linewidth]{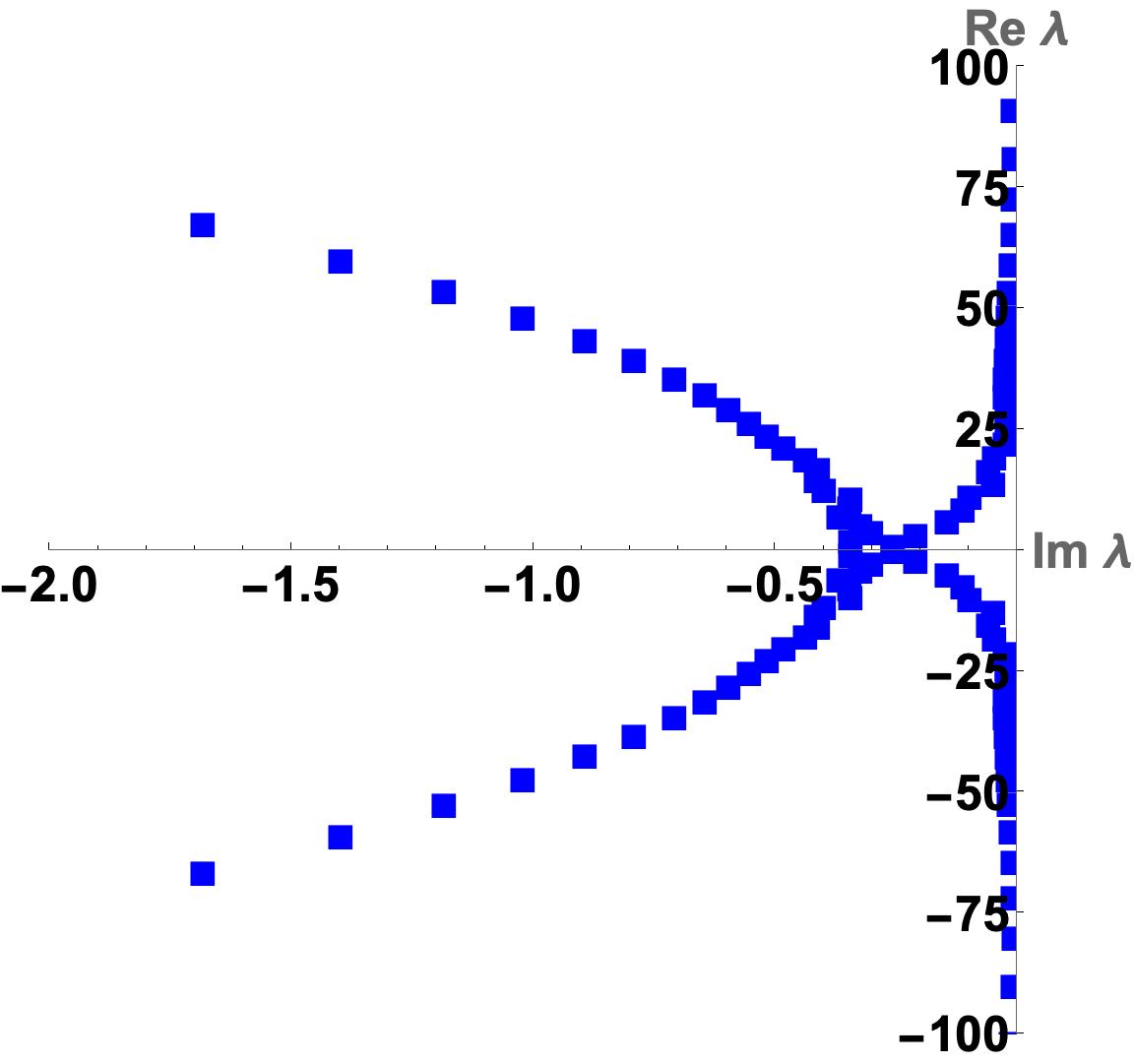}
    \end{minipage}~~
    \begin{minipage}{0.4\linewidth}
        \centering
        \small \textbf{(b)} $d_1 = 0$, $b_0 = b_1 = 1$\\[0.5em]
        \includegraphics[width=\linewidth]{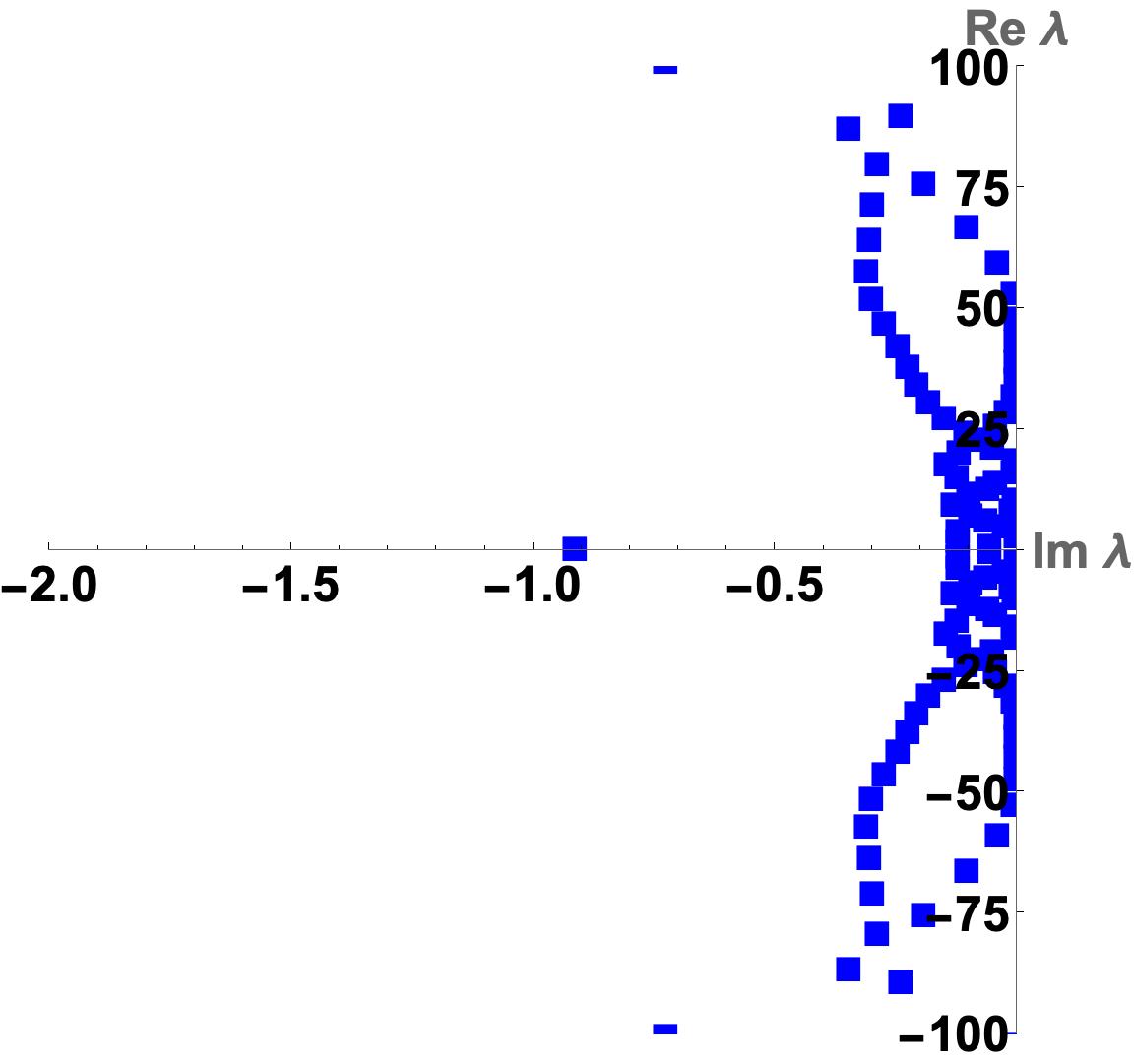}
    \end{minipage}

    \vspace{0.5em}

    \begin{minipage}{0.4\linewidth}
        \centering
        \small \textbf{(d)} $b_0 = b_1 = d_1 = 1$\\[0.5em]
        \includegraphics[width=\linewidth]{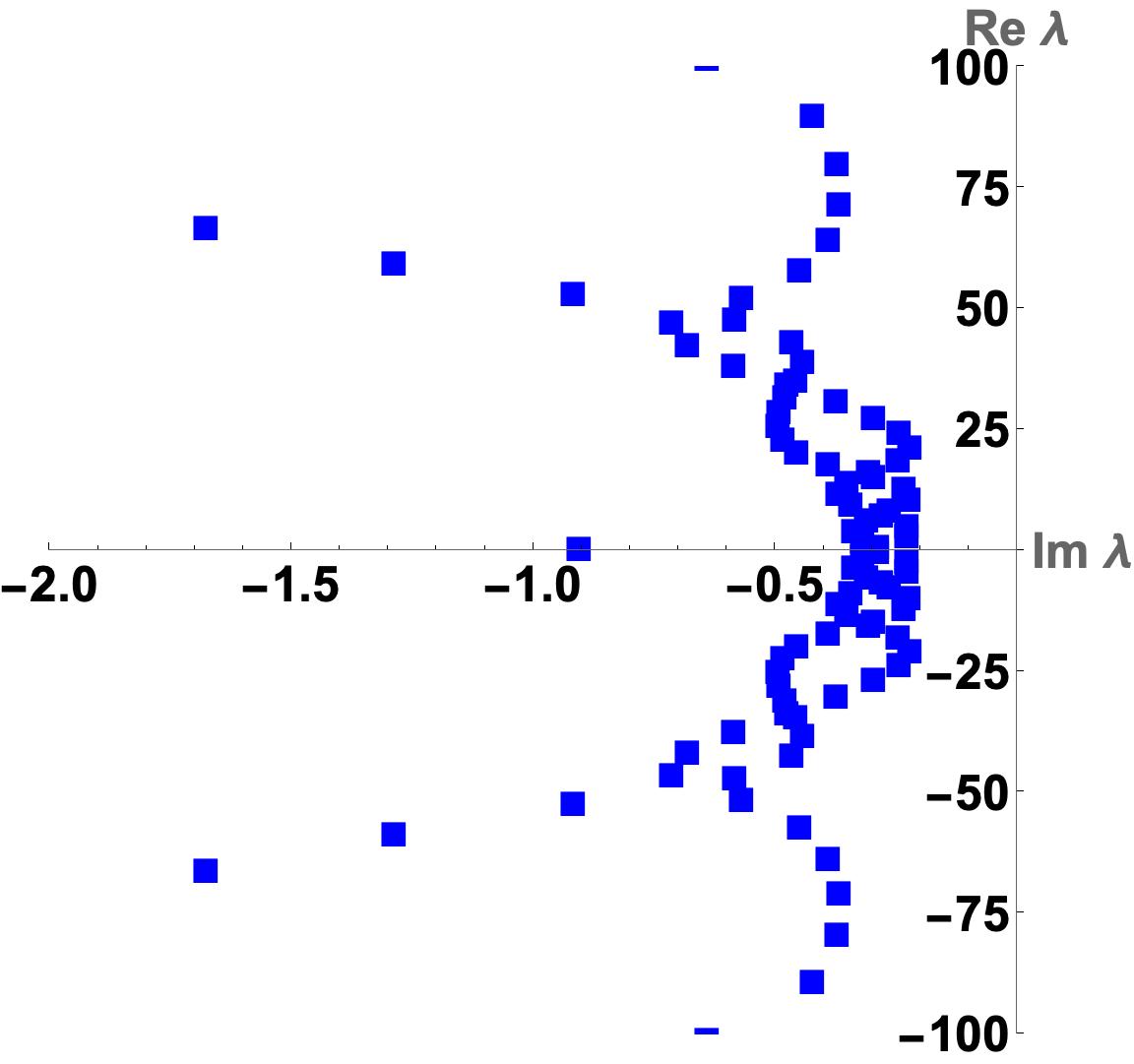}
    \end{minipage}~~
    \begin{minipage}{0.4\linewidth}
        \centering
         \small \textbf{(c)} $b_1 = 0$, $b_0 = d_1 = 1$\\[0.5em]
        \includegraphics[width=\linewidth]{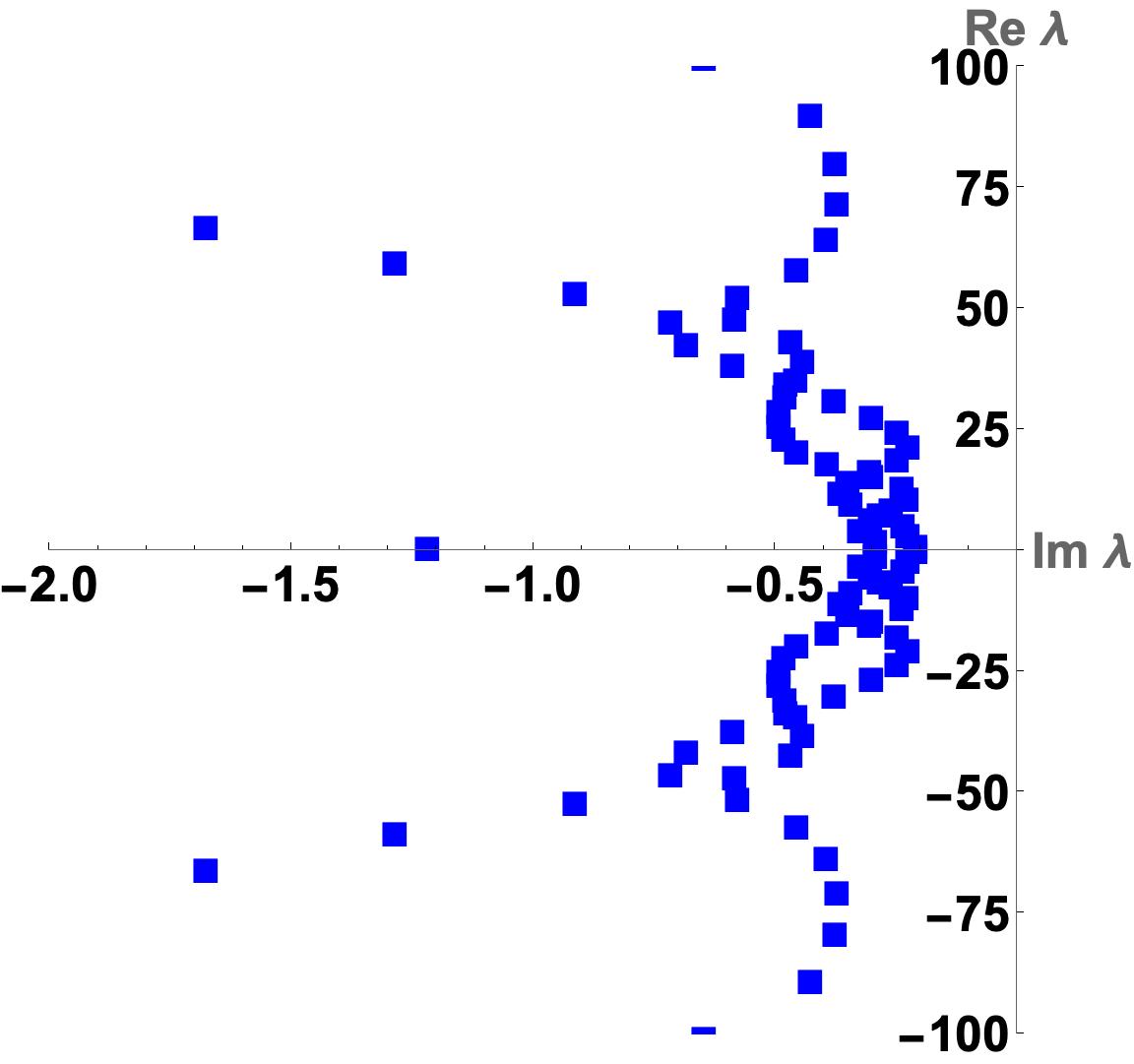}
    \end{minipage}

\caption{\footnotesize Eigenvalue distributions of the discrete system matrix $A_h$ under four damping configurations. Subfigures (a)--(d) illustrate how different combinations of the feedback gains $b_0$, $b_1$, and $d_1$ affect the spectral placement of eigenvalues and the resulting system stability.}
    \label{spectrals}
\end{figure}

\section{A Numerical Implementation}
\label{numerics}
To validate the theory, we run simulations under four damping setups, varying the feedback gains $b_0$, $b_1$, and $d_1$ to isolate the effects of higher-order, lower-order, and boundary damping. We use $\rho_1 = \sqrt{7}$, $\rho_2 = \pi$, $\alpha_1 = \sqrt{3}$, $\alpha_2 = 1$, $m = 0.6$, $b_0 = b_1 = d_1 = 1$, with $l_1 = 1$ and $l_2 = 2$, reflecting heterogeneity and highlighting damping influence. Simulations assess both the system’s spectral properties and its long-time behavior over a fixed horizon.

Following \cite{ozer2025uniformly}, let $N_1, N_2 \in \mathbb{N}$ with step sizes $h_1 = \frac{l_1}{N_1+1}$, $h_2 = \frac{l_2-l_1}{N_2+1}$ to discretize $[0,l_1]$ and $[l_1,l_2]$:
$0=x^1_0<\cdots<x^1_{N_1+1}=l_1$ and $l_1=x^2_0<\cdots<x^2_{N_2+1}=l_2.$
Let $x^i_{j+\frac{1}{2}} = (x^i_{j+1}+x^i_j)/2$ be midpoints. Approximating $w^i(x^i_j,t) \approx w^i_j(t)$ for $i=1,2$, $j=0,\ldots,N_i+1$, define
$
w^i_{j+\frac{1}{2}} := \frac{w^i_{j+1} + w^i_j}{2},$ $\delta_{x,h_i} w^i_{j+\frac{1}{2}} := \frac{w^i_{j+1} - w^i_j}{h_i},$ $ \delta_{x,h_i}^2 w^i_j := \frac{w^i_{j+1} - 2w^i_j + w^i_{j-1}}{h_i^2}.
$
First derivatives at midpoints improve accuracy via symmetric stencils \cite{ozer2025uniformly}.

A semi-discrete  approximation of \eqref{initialSys0}–\eqref{initialSys2} leads to
\begin{eqnarray} 
\left\{\begin{array}{ll}
\rho_i\frac{w^i_{j-\frac{1}{2}, tt} +w^i_{j+\frac{1}{2}, tt}}{2} - \alpha_i\delta_{x,h_i}^2 w^i_j = 0, ~~ j = 1, \dots, N_i \\
w^1_0 = 0, w^1_{N_1+1} = w^2_0(t) = z(t) \\
\alpha_1 \delta_{x,h_1} w^1_{N_1+\frac{1}{2}} + \rho_1 h_1 \frac{w^1_{N_1+1, tt} + w^1_{N_1, tt}}{4}- \alpha_2 \delta_{x,h_2}w^2_{\frac{1}{2}} \\  
 + \rho_2 h_2 \frac{w^1_{1, tt} + w^2_{0, tt}}{4} + m w^1_{N_1+1, tt} =  -b_1 w^1_{N_1+1, t} \\
 +b_0(\alpha_2 \delta_{x,h_2} w^2_{\frac{1}{2}, t} - \alpha_1 \delta_{x,h_1} w^1_{N_1+\frac{1}{2}, t}), \\
\alpha_2 \delta_{x,h_2} w^2_{N_2+\frac{1}{2}} + \rho_2 h_2 \frac{w^2_{N_2+1, tt} + w^2_{N_2, tt}}{4}  = -d_1 w^2_{N_2+1, t}\\
\left[w^i_j, w^i_{j,t}\right](0)=(w^i_0,w^i_1)(x^i_j,0), ~~i=1,2.
\end{array}
\right.
\end{eqnarray}

Let $\vec{\Phi} = [w^1_1,\dots,w^1_{N_1+1},w^2_1,\dots,w^2_{N_2+1}]^\top$ and $\vec{\Psi} := [\vec{\Phi},\vec{\Phi}_t]^\top$. The discretized system above can be written compactly as $\vec{\Psi}_t = A_h \vec{\Psi}$ where $A_h$ is the system matrix.

For illustration, we set $N_1 = N_2 = 30$. Figure~\ref{spectrals} shows the eigenvalue distributions of $A_h$ under four damping configurations.
{In subfigures (a) and (b), where $b_0 = 0$ or $d_1 = 0$, eigenvalues cluster near the imaginary axis, indicating polynomial decay \cite{Littman-Taylor,Avdonin-Edwards,Boughamda}. In contrast, (c) and (d) show uniform spectral separation when $b_0, d_1 > 0$, confirming exponential decay even without $b_1$. This highlights that higher-order and boundary damping alone suffice for exponential stability, refining earlier results that required lower-order damping.}


\begin{figure}[htb]
    \centering
    {\small \textbf{(a)} $b_0 = 0$, $b_1, d_1=1$}\\
    \includegraphics[width=0.42\linewidth]{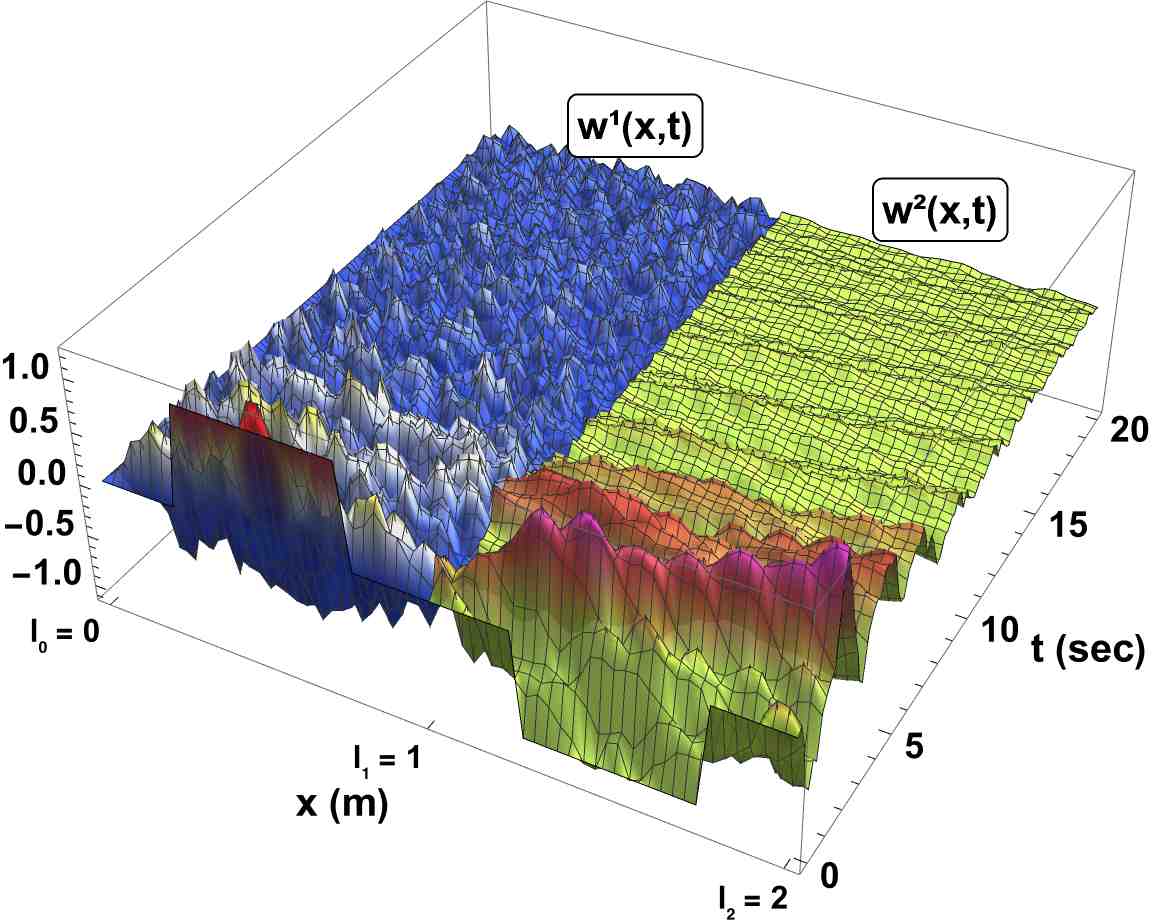}~~    
    \includegraphics[width=0.42\linewidth]{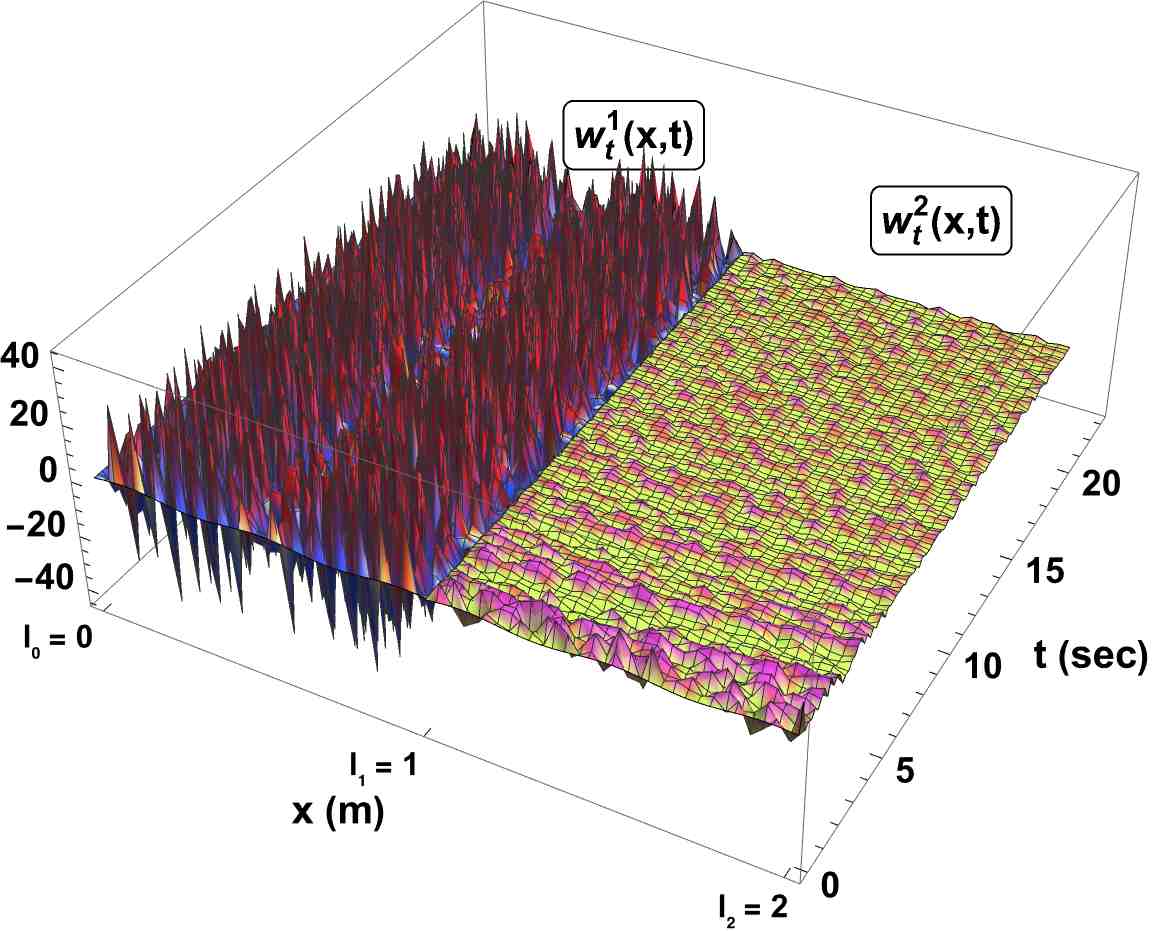}\\[.5em]
    {\small \textbf{(b)} $d_1 = 0$, $b_0, b_1=1$}\\
    \includegraphics[width=0.42\linewidth]{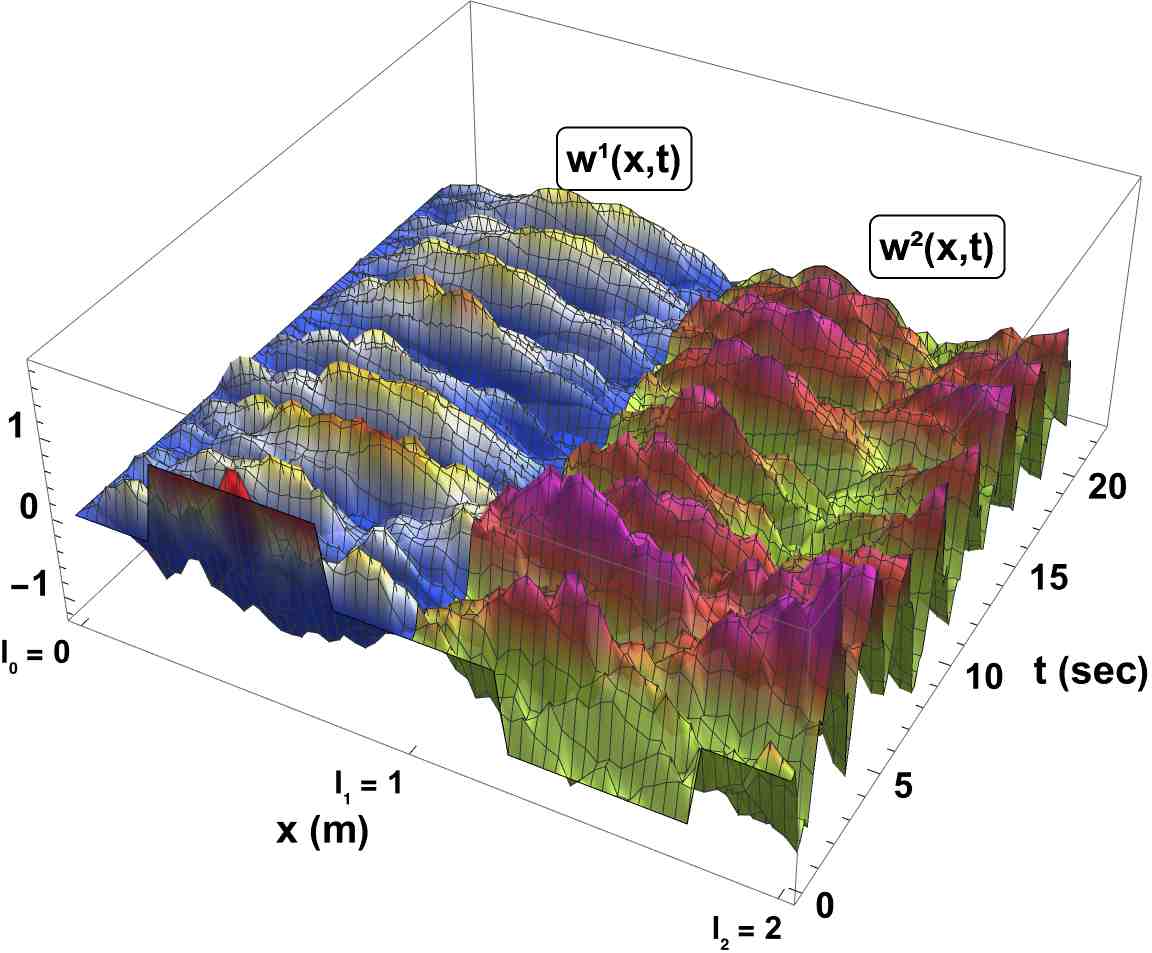}~~    
    \includegraphics[width=0.45\linewidth]{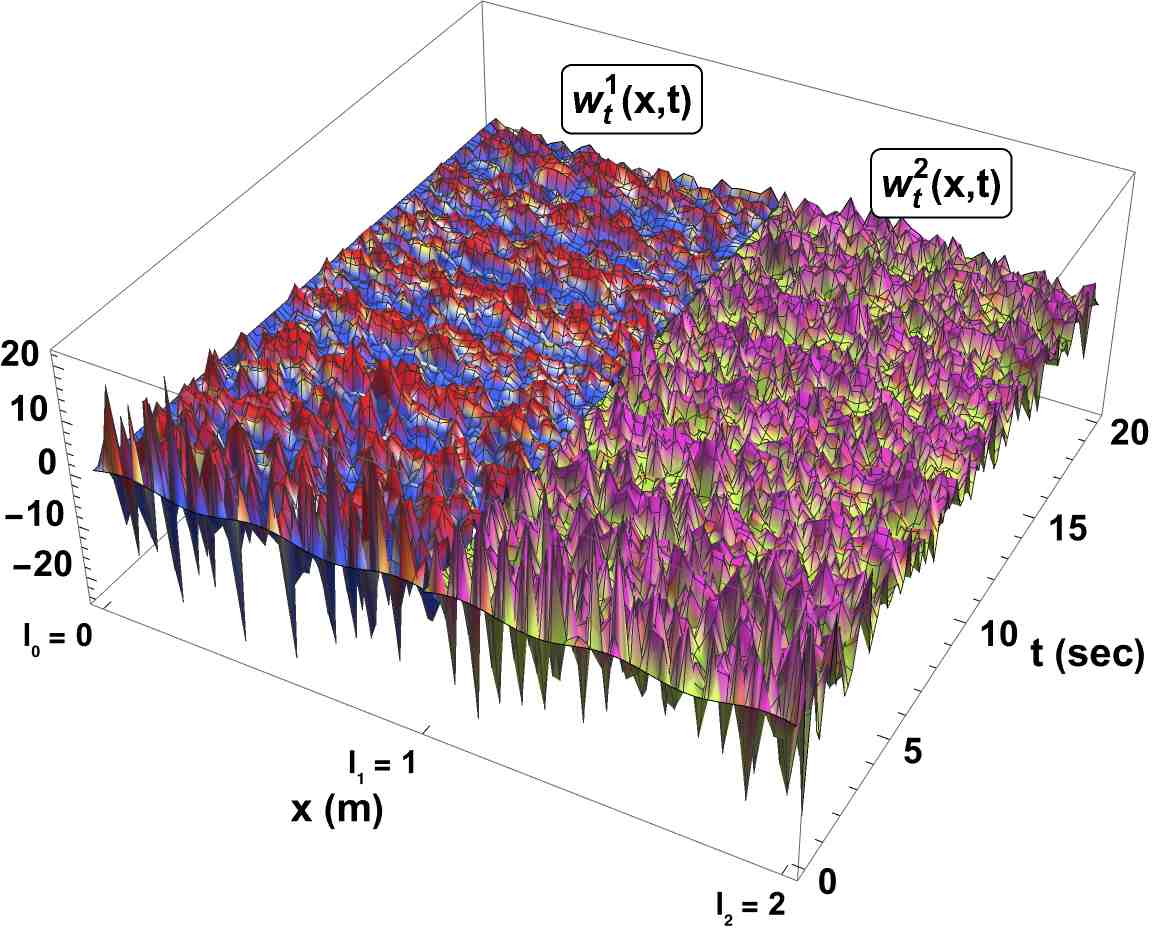}\\[.5em]
   \small \textbf{(c)} $b_0, b_1, d_1 = 1$\\
    \includegraphics[width=0.42\linewidth]{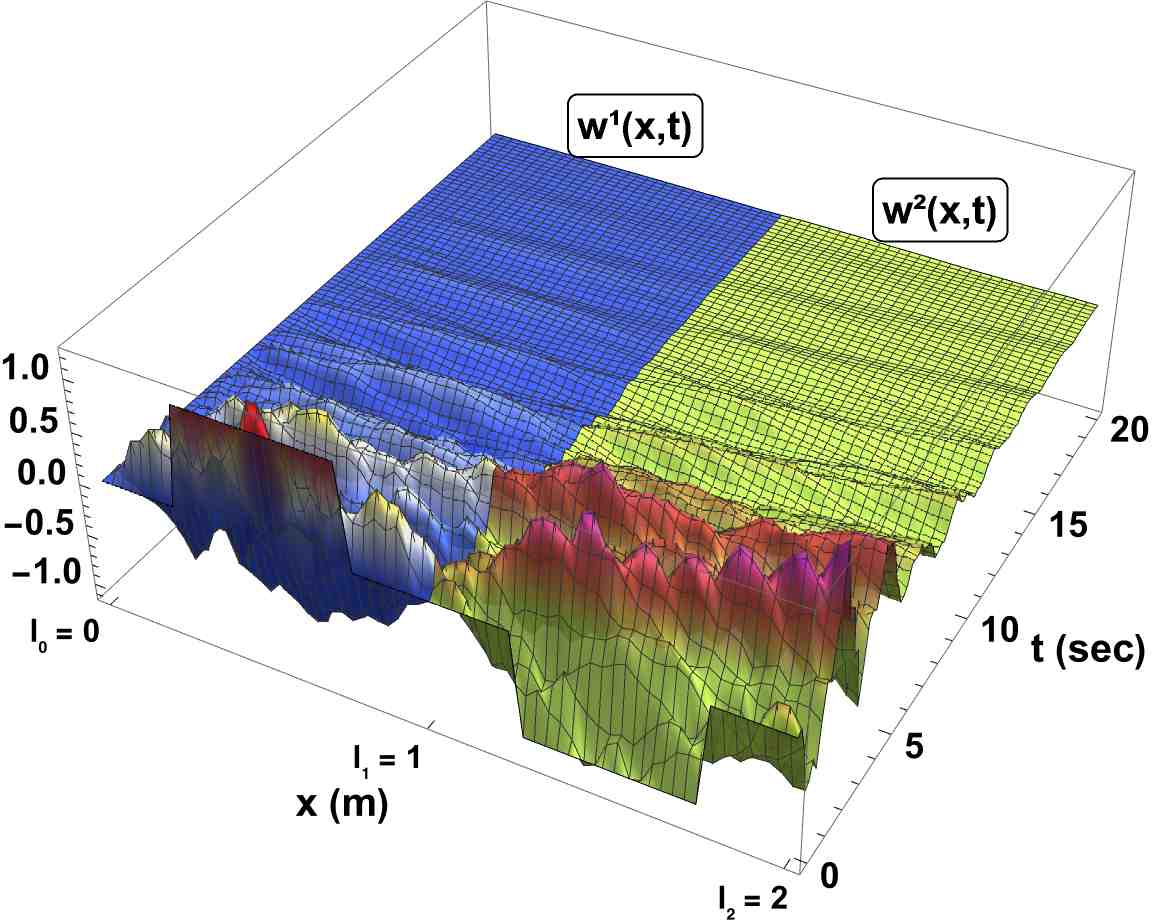}~~    
    \includegraphics[width=0.42\linewidth]{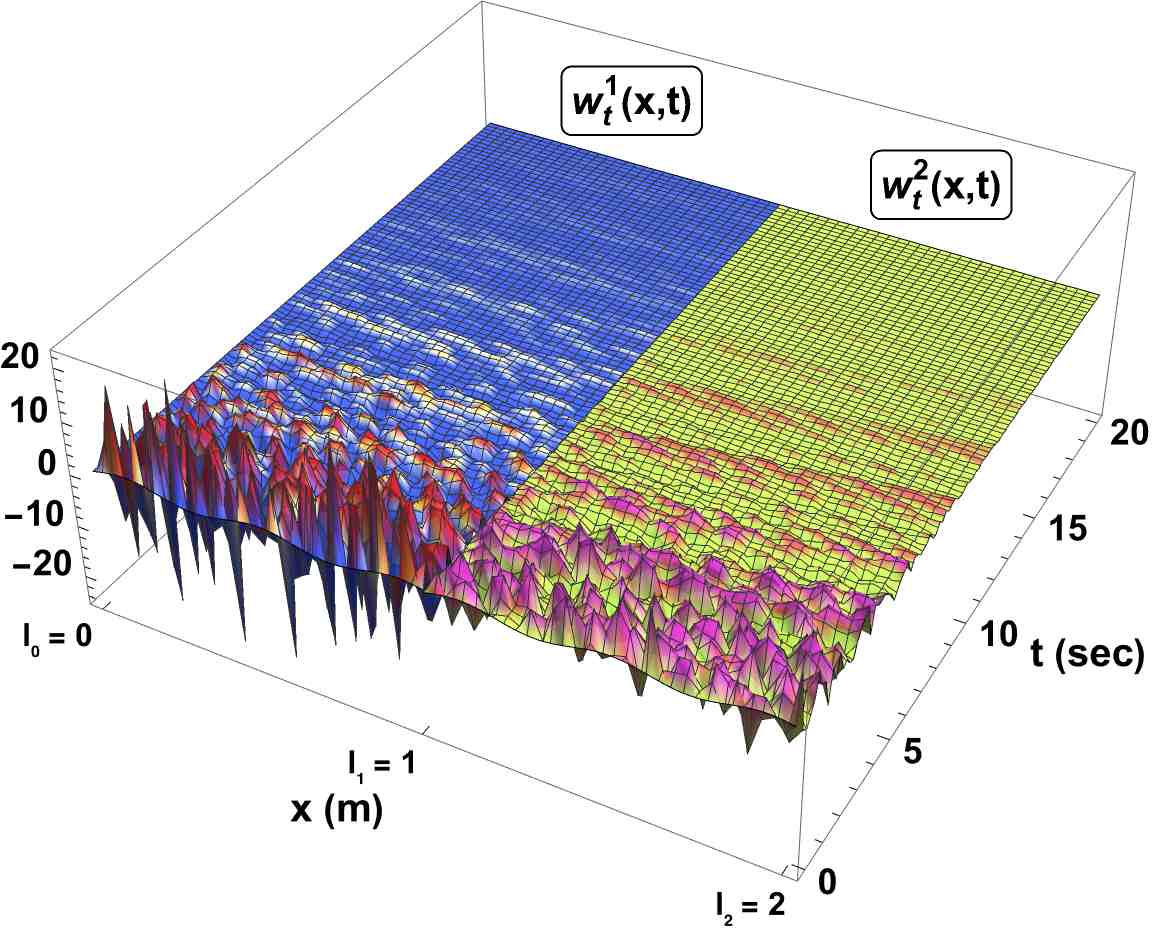}\\[.5em]
     \small  \textbf{(d)} $b_1 = 0$, $b_0, d_1=1$\\
    \includegraphics[width=0.42\linewidth]{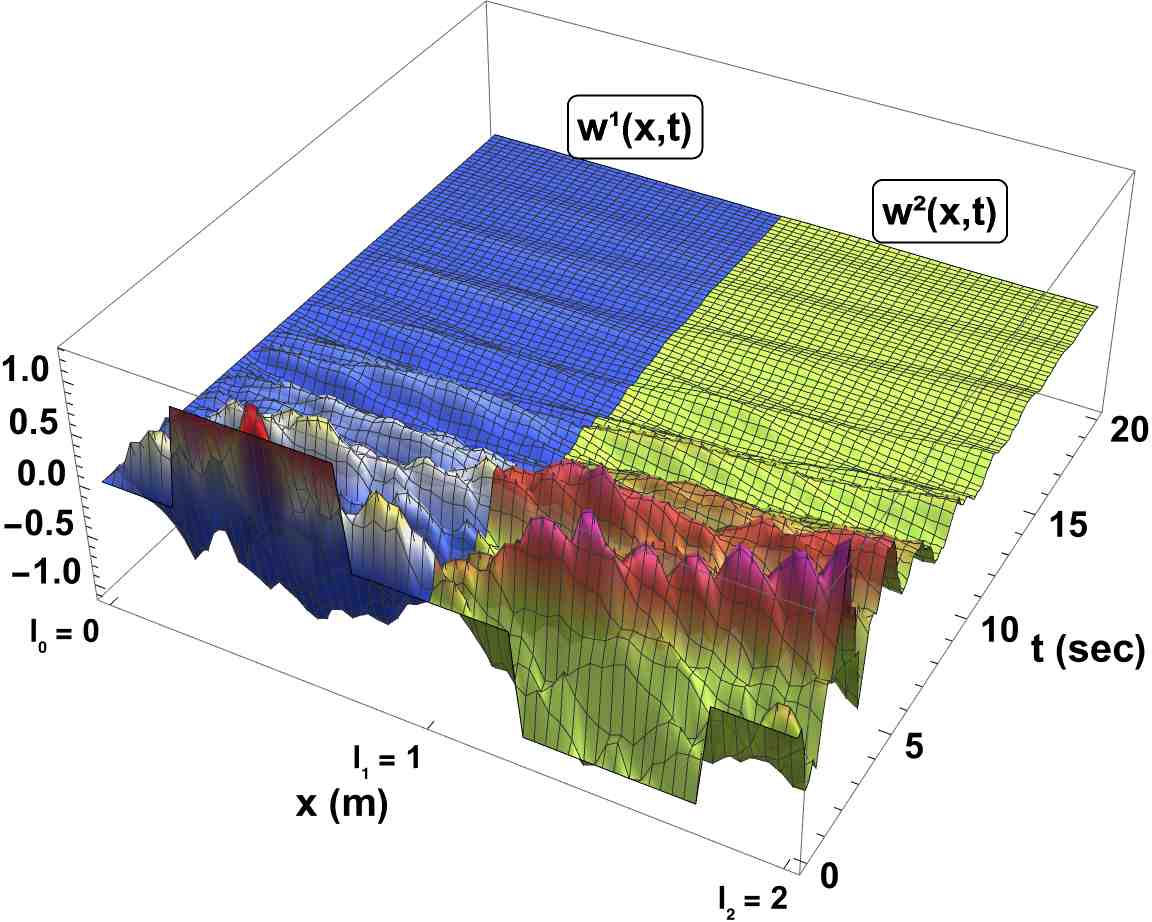}~~    
    \includegraphics[width=0.42\linewidth]{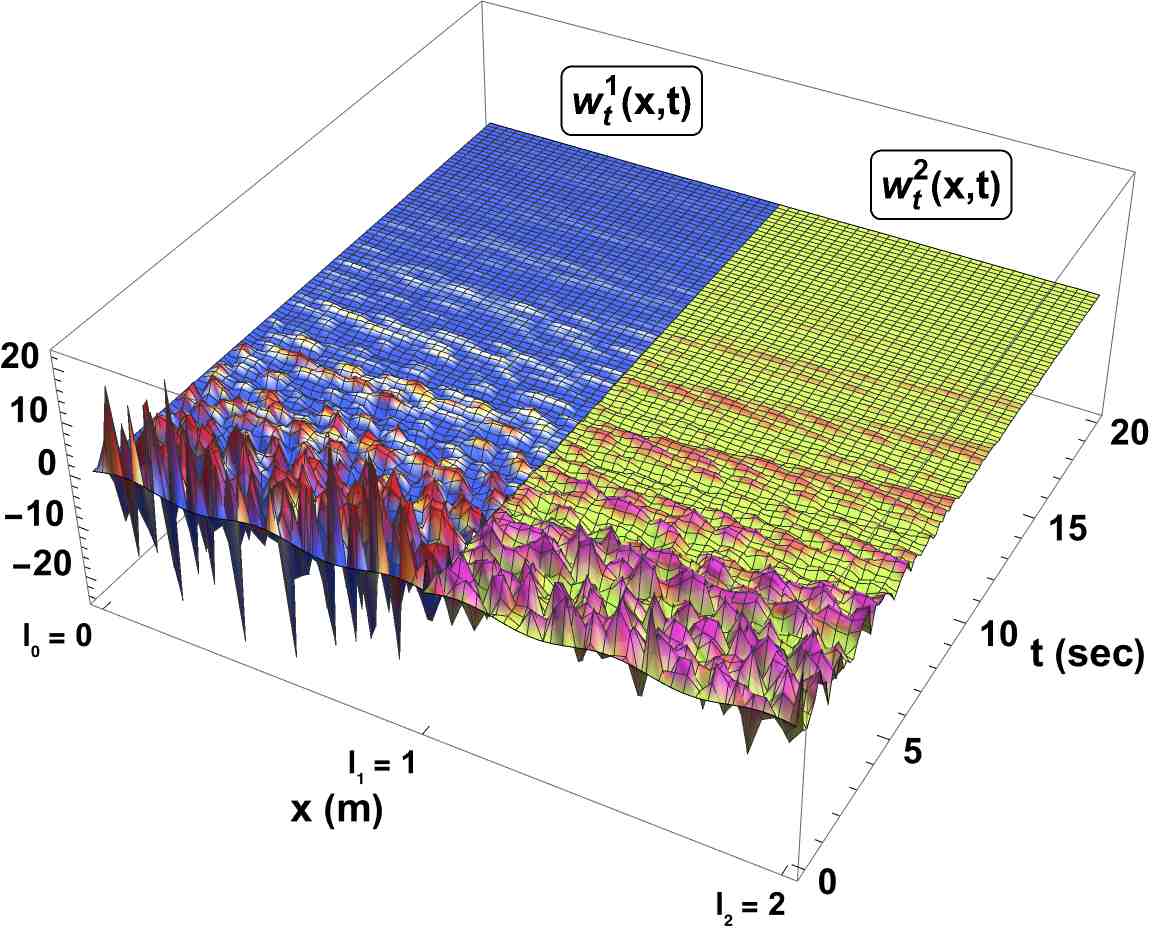}\\[.5em]
  \caption{\footnotesize Time evolution of the displacement fields $\{w^1(x,t), w^2(x,t)\}$ (left column) and their time derivatives $\{w^1_t(x,t), w^2_t(x,t)\}$ (right column) over a 20-second interval. Each row is labeled (a)--(d) and corresponds to a distinct damping configuration as indicated above.}
    \label{w1w2}
\end{figure}


\begin{figure}[htb]
    \centering
    \includegraphics[width=0.7\linewidth]{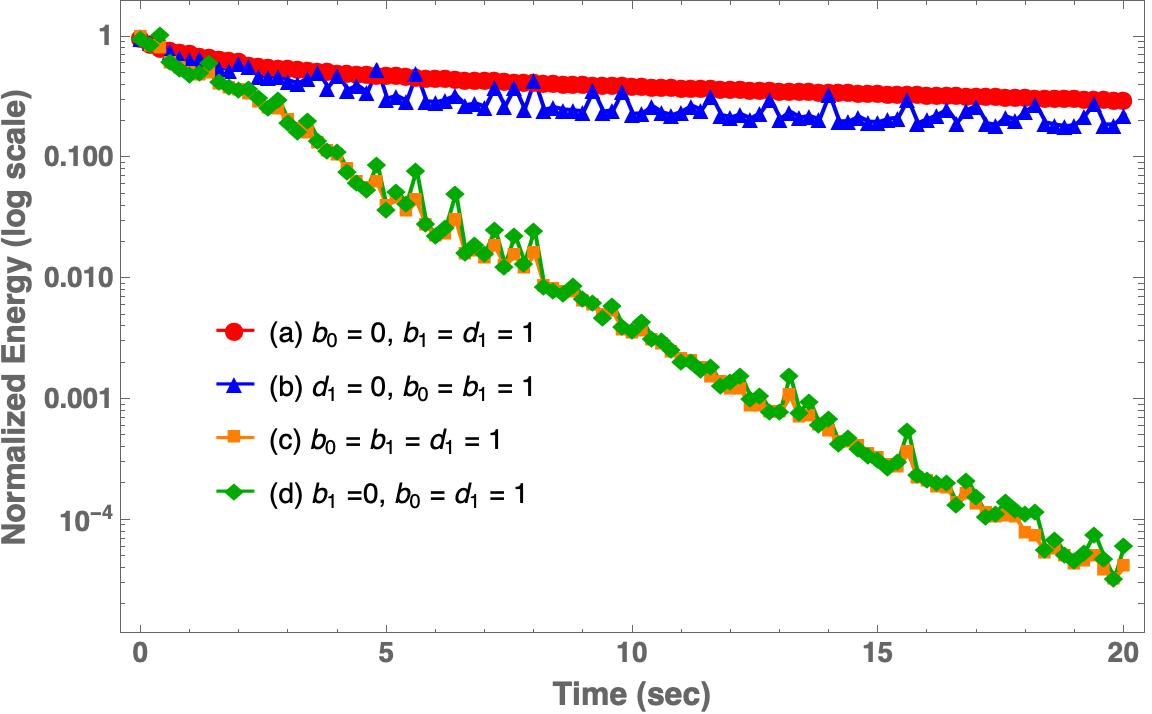} 
\caption{\footnotesize Normalized energy decay (log scale) over 20 seconds for damping configurations (a)--(d). Cases (a) and (b), lacking $b_0$ or $d_1$, show polynomial decay. Cases (c) and (d) confirm exponential decay, highlighting the role of higher-order and boundary damping.}
    \label{logenergy}
\end{figure}
\subsection{A Numerical Experiment}
To complement the spectral analysis in Fig.~\ref{spectrals}, we present time-domain simulations over a 20-second horizon under the same four damping configurations. These simulations visualize the decay of displacement and velocity fields and confirm the theoretical predictions. To assess robustness, we consider sinusoidal (continuous) and box-type (discontinuous) initial conditions. For $i = 1,2$, we set
$w^i_1(x) = \sin\left( \frac{4\pi x}{l_2} \right),$ and $
w^i_0(x) = (-1)^{i+1} \, \chi_{I_i}(x),$
where $I_i := \left( \dfrac{l_i - l_{i-1}}{2} - \dfrac{l_i - l_{i-1}}{4}, \, \dfrac{l_i - l_{i-1}}{2} + \dfrac{l_i - l_{i-1}}{4} \right)$ and $\chi_{I_i}$ is the characteristic function of $I_i$.
These allow us to examine the system's response to both smooth and localized initial data, highlighting the effectiveness of the damping configuration across varying excitation profiles.

Figure~\ref{w1w2} shows the time evolution of the displacement fields $w^1(x,t)$ and $w^2(x,t)$ over a 20-second interval for the four damping configurations (a)--(d). Figure~\ref{logenergy} displays the corresponding normalized energy decay.

{The numerical scheme, based on averaged Finite Differences, was recently developed for this coupled PDE–ODE model \cite{ozer2025uniformly} and is intentionally kept simple and explicit to ensure clarity and reproducibility; see also \cite{Guo3}. It reliably captures long-time energy behavior and highlights the distinct effects of each damping configuration.}


Configurations (a) and (b), with $b_0 = 0$ or $d_1 = 0$, lead to slow, non-exponential decay, confirming that {isolated or lower-order damping is inadequate}, consistent with Figure~\ref{spectrals}. Configuration (c), with all damping active, achieves rapid exponential decay (Theorem~\ref{majres}). Configuration (d), with $b_1 = 0$, still yields exponential decay, showing that $b_0$ and $d_1$ suffice (Theorem~\ref{majres2}). These results show effective stabilization is driven by {boundary and higher-order damping}.

\begin{remark}
The simulations use $\epsilon_1 = 0.07087$, $\epsilon_2 = 0.12073$ (within Theorem~\ref{majres0} bounds), and $\delta = 1.53515$, yielding $T = \max(T_1, T_2, T_3, T_4) = 70.22$. Yet, stabilization occurs within 20 seconds, confirming the control’s efficiency.
\end{remark}

\section{Conclusions \& Future Work}
We prove exponential stability for two strings joined by a dynamic interior mass. Our analysis shows that stability holds without the lower-order damping $b_1$, allowing a simpler, more robust design. The result is unconditional, independent of wave speeds or interface location.



A companion study \cite{Akil} classifies decay rates under six damping scenarios using operator-theoretic tools, covering exponential and polynomial regimes. Future work includes extensions to $k \geq 2$ serially-connected wave equations with multiple interior masses, {under partial damping, discontinuous joints \cite{Walterman2025}, and/or sparse actuators}. Other challenges include dynamic boundaries, time-delay feedback, and control-matched disturbances, as in \cite{Mei, nicaise2011interior}.


\bibliographystyle{IEEEtran}
\bibliography{references}
\end{document}